\documentclass[12pt]{article}

\usepackage[T2A]{fontenc}
\usepackage[utf8]{inputenc}
\usepackage[english]{babel}
\usepackage[tbtags]{amsmath}
\usepackage{amsfonts,amssymb,mathrsfs,amscd}

\usepackage{hyperref}

\newtheorem{theorem}{Theorem}
\newtheorem{lemma}{Lemma}

\newtheorem{proposition}{Proposition}
\newtheorem{corollary}{Corollary}
\newtheorem{remark}{Remark}

\def\E{\hskip.15ex\mathsf{E}\hskip.10ex}
\def\P{\mathsf{P}}

\title{On efficient estimates of the rate of convergence for Markov chains}
\author{A.Yu. Veretennikov\footnote{Institute for information transmission problems RAS (Kharkevich Institute); \\
email: ayv@iitp.ru 
} }

\date{}

\begin{document}

\maketitle



\begin{abstract}
The paper presents efficient approaches for evaluating convergence rate in total variation for finite and general 
Markov chains. The motivation for studying convergence rate in this metric is its usefulness in various limit theorems (see, for example, \cite{Dobr56}).
For homogeneous Markov chains the goal is to compare several different methods: (1) the second eigenvalue for the transition matrix method (the method no. 1), (2) the method based on Markov -- Dobrushin's  ergodic coefficient, and the new spectral method developed in \cite{VV_Dokl, VerVer}, as well as modifications of they both by iterations (the ``other methods''). We answer the question whether or not the ``other methods''  may provide the optimal or close to optimal convergence rate in the case of homogeneous Markov chains. The answer turns out to be positive for appropriate modifications of both ``other methods''. The analogues of these ``other methods'' for the non-homogeneous Markov chains are also presented.

The work is theoretical.  However, the methods of computing efficient bounds of convergence rates may be in demand in various applied areas. 

\end{abstract}

\section{Introduction}
For  Markov chains (MC) the new spectral approach for evaluating the rate of convergence in total variation uses a so called Markovian coupling. Certain sub-stochastic matrices or operators\footnote{The structure of the operators \(V^{(m)}\) will be explained in what follows, see (\ref{V}).}, \(V\) and \(V^{(m)}\) are introduced with the help of this coupling construction, and the rate of convergence under investigation is estimated by their spectral radii \(r(V), r(V^{(m)})\). 
The norm \(\|V\|\) coincides with the Markov -- Dobrushin ergodic coefficients (MD in what follows, see the definition in (\ref{MD}) in the next section), although, the formal definitions may seem quite different. The advantage of using spectral radii is due to the well-known fact that in all cases \(r(V^{(m)}) \le \|V^{(m)}\|\), so that the new proposed bounds are usually better. 

These several methods will be compared, and the analogues of the estimates based on the operators \(V\) and \(V^{(m)}\)  will be presented for non-homogeneous Markov chains where the estimates based on the spectral radii are not applicable. 
By the suggestion of the anonymous referee, the rest of this introduction is moved to the conclusion remarks section.

The paper consists of four  sections: 1 -- this short introduction, 2 -- homogeneous Markov coupling construction and sub-Markov operators, 3 --  non-homogeneous case; 4 --  conclusive remarks. The main results are the theorems \ref{thm1} --  \ref{thm3}. 


\section{Homogeneous coupling and operators $V, V^{(m)}$}
In this section we recall briefly the main points of the presentation from \cite{VerVer} with certain changes required for the new setting. 
Also, a short presentation of \cite{VerVer} may be found in \cite{VV_Dokl}. Unlike \cite{VV_Dokl, VerVer}, here the emphasis will be on the iterative $m$-step coupling, which was only  briefly outlined in \cite{VerVer}. 

\medskip

\subsection{Markov -- Dobrushin's and Gantmacher's methods}\label{sec:homo} 
We consider firstly a homogeneous Markov process (MP) in discrete time $(X_n, \, n\ge 0)$ in a general non-empty state space $S$ with a topology and with a Borel sigma-algebra ${\cal S}=\sigma(S)$; as usual in Markov processes, any state $\{x\}$ belongs to $\sigma(S)$ by assumption.  
If  the state space $S$ is finite, then $|S|$ denotes the number of its elements, and $\cal P$ stands for the transition matrix $\left(p_{ij}\right)_{1\le i,j \le |S|}$ of the process. The latter notation may also be applied in the case where $S$ is countable. Throughout the paper, the chain is assumed irreducible and aperiodic (also called primitive). By $\mu^x_n(dy)$ we understand the marginal measure of the process at time $n$ given the initial state $x$; $\P_x(X_n\in dy)$ is an equivalent notation; $\mu_{inv}$ will stand for the invariant measure of the process. Also, $\P_x(m,dy)$ denotes the $m$-step transition kernel of the process ``from $x$ to $dy$''. The Markov--Dobrushin constant (complementary to Markov -- Dobrushin's ergodic coefficient, see what follows) is defined by the formula 
\begin{equation}\label{MD}
\alpha:= 1-\frac12\,\sup_{x,x'}\|\mu^x_1 - \mu^{x'}_1\|_{TV} = \inf_{x,x'} \int \left(\frac{\P_{x'}(1,dy)}{\P_x(1,dy)}\wedge 1 \right)\P_x(1,dy). 
\end{equation}

Recall that the following well-known inequality in the ergodic theorem -- which was established by A.A. Markov himself in the case where $S$ is finite and all transition probabilities are positive --  assuming that $\alpha>0$,

\begin{equation}\label{bound1}
\sup_{x}\sup_{A\in S} |\mu^x_n(A) - \mu_{inv}(A)| \le  (1-\alpha)^{n}, \quad n\ge 0. 
\end{equation}
Notice that here $\displaystyle \frac{\P_{x'}(1,dy)}{\P_x(1,dy)}$ in the right hand side of (\ref{MD}) is understood in the sense of the density of the absolute continuous component of the numerator with respect to the denominator. In other words, the two measures are not necessarily have to be equivalent, but may possess singular components with respect to each other. The same relates to the expression $\displaystyle \frac{\P_{x'}(m,dy)}{\P_x(m,dy)}$ in what follows. 
Denote also, 
$$
\delta:=1-\alpha.
$$
This $\delta$ is Markov -- Dobrushin's ergodic coefficient.

The inequality (\ref{bound1}) may be extended in the  proposition below as follows (see \cite[proposition 5]{VerVer}). Recall that $X_n$ is a Markov chain in a general state space, and only the references to Gantmacher's result will relate to the case of a finite $S$.

In the next proposition the $m$-step analogue of the coefficient $\alpha$ is used:
\begin{equation}\label{MDm}
 \alpha^{(m)} := \inf_{x,x'} \int \left(\frac{\P_{x'}(m,dy)}{\P_x(m,dy)}\wedge 1 \right)\P_x(m,dy), \quad m\ge 1. 
\end{equation}
Recall that the process is irreducible and aperiodic.

\begin{proposition}[\cite{VerVer}, proposition 5]
\label{thm1} 
Let the homogeneous Markov chain $(X_n)$ be irreducible and aperiodic, and let  $m\ge 1$ be such that $ \alpha^{(m)} >0$. 
Then the process $(X_n)$ is ergodic, i.e., there exists  a limiting invariant probability measure $\mu_{inv}$, which is stationary and for every $n$, 
and for any $m\ge 1$
\begin{equation}\label{pro1-e1}
\sup_{x}\sup_{A\in S} |\mu^x_n(A) - \mu_{inv}(A)| \le  (1-\alpha^{(m)})^{[n/m]}. 
\end{equation}
\end{proposition}
Further, denote 
$$
\delta^{(m)} := 1 - \alpha^{(m)}.
$$
This $\delta^{(m)}$ is iterated  Markov -- Dobrushin's (MD) $m$-step ergodic coefficient. 
Notice elementary inequalities,  
$$
\delta^{(m)} \le \delta^m, 
$$
for any $m\ge 1$, 
due to the Markov property. In particular, for any two integers $m_1,m_2\ge 1$, if we let $m_3:= m_1m_2$, then for any $k\ge 1$ we have, 
\begin{equation}\label{dm123}
(\delta^{(m_3)})^{k}\le (\delta^{(m_i)})^{(km_3/m_i)} , \quad i=1,2.
\end{equation}
Therefore, up to some bounded multipliers, the larger is $m$, the more precise -- at least, asymptotically -- is the bound in the exponential bound (\ref{pro1-e1}). 

It is important to notice that there are natural examples where 
$\delta=0$, in which case this coefficient $\delta$ does not provide any efficient bound for the rate of convergence, while $\delta^{(2)}>0$, which does provide such a bound; examples may be found in \cite{NV1}.

Further, for any Borel measurable $A$, the function  $\P_x(X_1\in A)$ is Borel measurable with respect to $x$; this is a standard requirement in Markov processes \cite{Dynkin}. Such measurability with respect to the pair $x,x'$ will be also valid for the measure $\Lambda^{(m)}_{x,x'}$ defined below due to the linearity. 
For two fixed states $x, x'$ and for any $m\ge 1$, denote 
\begin{equation*}
\Lambda^{(m)}_{x,x'}(dz) := \P_x(m,dz) + \P_{x'}(m,dz).
\end{equation*} 
Clearly,  $\Lambda^{(m)}_{x,x'}(dz) = \Lambda^{(m)}_{x',x}(dz)$.
Then the definition of the coefficient $\alpha^{(m)}$ may be rewritten in the form, 
$$
 \alpha^{(m)} = \inf_{x,x'} \int \Big(\frac{\P_{x'}(m,dy)}{\Lambda^{(m)}_{x,x'}(dy)}\wedge \frac{\P_{x}(m,dy)}{\Lambda^{(m)}_{x,x'}(dy)} \Big)\Lambda^{(m)}_{x,x'} (dy), 
$$
which is clearly  symmetric with respect to $x$ and $x'$. In fact, a similar formula 
gives the same result for any other measure with respect to which both $\P_{x'}(m,dy)$ and $\P_{x}(m,dy)$ are absolutely continuous, because if $d\Lambda^{(m)}_{x,x'} <\!\!< d\tilde\Lambda^{(m)}_{x,x'}$, 
then 
\begin{align*}
&\displaystyle  \int \!\Big(\frac{\P_{m,x'}(dy)}{\Lambda^{(m)}_{x,x'}(dy)}\wedge \frac{\P_{x}(m,dy)}{\Lambda^{(m)}_{x,x'}(dy)} \Big)\Lambda^{(m)}_{x,x'} (dy) 
\!= \!\int\! \Big(\frac{\P_{x'}(m,dy)}{\tilde\Lambda^{(m)}_{x,x'}(dy)}\wedge \frac{\P_{x}(m,dy)}{\tilde\Lambda^{(m)}_{x,x'}(dy)} \Big)\tilde\Lambda^{(m)}_{x,x'} (dy). 
\end{align*}

\medskip

Here is one of the important notions in the following presentation: for any $m\ge 1$, denote 
\[
\alpha^{(m)}(x,x') : = \int \left(\frac{\P_{x'}(m,dy)}{\P_x(m,dy)}\wedge 1 \right)\P_x(m,dy).
\]
Clearly, for any $x,x'\in S$, 
\begin{equation*}
\alpha^{(m)}(x,x') = \alpha^{(m)}(x',x) \ge \alpha^{(m)}. 
\end{equation*}
This function $\alpha^{(m)}(x,x')$ will be used to construct sub-Markovian operators $V^{(m)}$ and $\hat V^{(m)}$. These operators play a key role in the   extensions of the bounds (\ref{bound1}) and  (\ref{pro1-e1}), see theorems \ref{thm2} and \ref{thm3} in what follows. 

\medskip

In the meanwhile, let us now recall a simplified version of the Gantmacher formula \cite[formula (13.96)]{Gant},  based on the second eigenvalue and other spectral properties of the transition matrix in the case of a finite $S$. This formula in \cite{Gant}, actually, provides a precise, non-asymptotical result\footnote{It may be understood from the comments in \cite[ch. XIII, \S 7]{Gant} that the motivation for the author was the  question  posed by Kolmogorov. } in terms of the full  spectral decomposition of the transition matrix. We only state here its corollary. The notation $\lambda_2({\cal P})$ stands for the second eigenvalue of the transition matrix~${\cal P}$.

\begin{proposition}[\cite{Gant}, corollary of the formula (13.96) in ch. XIII, \S 7]\label{pro2}
If $S$ is finite and  the transition matrix ${\cal P}$ is primitive, then the convergence to the invariant measure satisfies the bound 
\begin{equation}\label{Gantf}
C |\lambda_2({\cal P})|^n-o(|\lambda_2({\cal P})|^n) \le \|\mu^{x}_n - \mu_{inv}\|_{TV}\le C |\lambda_2({\cal P})|^n + o(|\lambda_2({\cal P})|^n),
\end{equation}
for $n$ large enough, uniformly with respect to the initial state $x$; here $C=C({\cal P})$ admits a representation in terms of the algebraic structure of the transition matrix $\cal P$. 

\end{proposition}

As we shall see shortly, the formula (\ref{Gantf})  implicitly prompts that by using the MD method, or the intermediate bound by the spectral radii of the sub-Markov operators $V^{(m)}$, or $\hat V^{(m)}$ (see what follows), it is also possible to achieve similar bounds for the rate of convergence if $m$ is chosen large enough. 
  
\medskip  

Yet, it should be recalled that  for non-homogeneous MC this best recipe based on the second eigenvalue of the transition matrix is not working, just because all transition matrices may differ at different times, hence, be non-commutative and have no joint spectral radius. On the other hand, the MD, iterated MD (i.e., based on the coefficient $\alpha^{(m)}$ with $m>1$), and $(V^{(m)})$ methods (see what follows) work well and allow to improve the ``simple'' MD rate of convergence bound. They  also work well in general state spaces where the applicability of the Gantmacher method is unclear. 

In the next result the state space $S$ may be general, not necessarily finite. However, in the case of a finite $S$, the role of $\lambda$ in the condition (\ref{thm1-e1}) in the theorem \ref{thm1} plays $\lambda = |\lambda_2({\cal P})|$.

\begin{theorem}\label{thm1} 
Let the homogeneous Markov chain $(X_n)$ be irreducible and aperiodic, and let  $m\ge 1$ be such\footnote{It is known that, indeed, for irreducible and aperiodic chains $ \alpha^{(m)} \!>\!0$ starting from some $m_0$.} that $ \alpha^{(m)} >0$. 
If for some $C>0$ and $\lambda<1$  
there exists $m_0$ such that 
\begin{equation}\label{thm1-e1}
\sup_{x_1,x_2} \|\mu^{x_1}_m - \mu^{x_2}_m\|_{TV}\le C \lambda^m, \quad \forall \, m\ge m_0, 
\end{equation} 
then $\forall \, \varepsilon>0$ there exists  $m_1\ge m_0$ 
such that 
(let $\lambda_\varepsilon = (1+\varepsilon)\lambda$)
\begin{equation}\label{thm1-e2}
(1-\alpha^{(m_1)}) 
\le 
(1+ \varepsilon/2)^{m_1} \lambda^{m_1} = \lambda_{\varepsilon/2}^{m_1}.
\end{equation}
and for any $n\ge m_1$, 
\begin{align}\label{thm1-e3}
&\sup_{x_1,x_2}\|\mu^{x_1}_n - \mu^{x_2}_n\|_{TV} 
\le \lambda_{\varepsilon}^n, \quad \forall \, n\ge m_1,
\end{align}
and also, 
\begin{equation}\label{thm1-e4}
\sup_{x_1} \|\mu^{x_1}_n - \mu_{inv}\|_{TV} \le  
\lambda_\varepsilon^{n}, \quad n\ge m_1, 
\end{equation}
where $\mu_{inv}$ is the unique limiting invariant measure. 
\end{theorem}
{\bf Proof} is nearly a tautology, and is typed here for completeness.  
Of course, it only makes sense to consider $\varepsilon$ such that $(1+\varepsilon)\lambda<1$.
The first claim is elementary. Let us choose $m'_1\ge m_0: \, (1+\varepsilon)^{m'_1}\ge C$, that is, $m'_1 \ln (1+\varepsilon) \ge (\ln C) \vee  m_0$. The value $m_1 \ge m'_1$ will be fixed a bit later. In any case, we have, 
$$
\sup_{x_1,x_2} \|\mu^{x_1}_{m'_1} - \mu^{x_2}_{m'_1}\|_{TV}\le (1+\varepsilon/2)^{m'_1}\lambda^{m'_1}, 
$$
and, clearly, the same inequality holds with any $m_1\ge m'_1$ in place of $m'_1$.

\medskip

Recall the notation $\lambda_\varepsilon = (1+\varepsilon)\lambda$. 
Then we get, 
\begin{equation}\label{thm1pr-e1}
(1-\alpha^{(m_1)}) 
\,= \sup_{x_1,x_2} \|\mu^{x_1}_{m_1} - \mu^{x_2}_{m_1}\|_{TV}\, \le 
(1+ \varepsilon/2)^{m_1} \lambda^{m_1} = \lambda_{\varepsilon/2}^{m_1},
\end{equation}
as required.

The proof of the remaining claims is as follows. 
We have,  
\begin{equation}\label{thm1pr-e2}
\sup_{x_1,x_2} \|\mu^{x_1}_{n} - \mu^{x_2}_{n}\|_{TV}\, 
\le (1-\alpha^{(m'_1)})^{[n/m'_1]}  \le 
\lambda_{\varepsilon/2}^{m'_1 [n/m'_1]} = (1+\varepsilon/2)^{m'_1 [n/m'_1]} \lambda^{m'_1 [n/m'_1]}.
\end{equation}
At this stage, let us choose $m_1$ so that, in addition to the requirement $m_1\ge m'_1$, it also satisfies 
$$
(1+\varepsilon/2)^{m'_1 [n/m'_1]} \lambda^{m'_1 [n/m'_1]} \le 
(1+\varepsilon)^{n} \lambda^{n}, \quad \forall \; n\ge m_1.
$$
Denote 
$$
C_{m'_1,\varepsilon} := (1+\varepsilon/2)^{- m'_1} \lambda^{-m'_1}.
$$
Then what we need is the inequality 
\begin{align*}
&(1+\varepsilon/2)^{m'_1 [n/m'_1]} \lambda^{m'_1 [n/m'_1]} 
\le C_{m'_1,\varepsilon}(1+\varepsilon/2)^{n} \lambda^{n}.  
\end{align*}
So, it suffices that $\forall \, n\ge m_1$,
$$
C_{m'_1,\varepsilon}(1+\varepsilon/2)^{n} \le (1+\varepsilon)^{n},
$$
or, equivalently, 
\begin{align*}
n \,\ln \frac{1+\varepsilon}{1+\varepsilon/2} \ge \ln C_{m'_1,\varepsilon}.
\end{align*}
Let us choose
$$
m_1 \ge  \Big\lceil \Big(\ln \frac{1+\varepsilon}{1+\varepsilon/2}\Big)^{-1} \ln C_{m'_1,\varepsilon}  \Big\rceil.
$$
Then, clearly, it follows from (\ref{thm1pr-e2}) that for any $n\ge m_1$ 
\begin{equation*}
\sup_{x,x'} \|\mu^x_{n} - \mu^{x'}_{n}\|_{TV}\, 
\le (1+\varepsilon)^{n} \lambda^{n} = \lambda_\varepsilon^n,
\end{equation*}
as required for the bound (\ref{thm1-e3}).

\smallskip

Finally, to show the estimate (\ref{thm1-e4}),  we have with the same $m_1$ and $n\ge m_1$, due to the triangle inequality and Chapman -- Kolmogorov formulae, 
\begin{align*}
&\sup_{x_1} \|\mu^{x_1}_n - \mu_{inv}\|_{TV} = 
\sup_{x_1} \|\mu^{x_1}_n - \int \mu^{x'}_n \mu_{inv}(dx')\|_{TV}
 \\\\
& = 
\sup_{x_1}  \| \int (\mu^{x_1}_n -\mu^{x'}_n) \mu_{inv}(dx')\|_{TV}
\le \sup_{x_1} \int \| (\mu^{x_1}_n -\mu^{x'}_n)\|_{TV} \mu_{inv}(dx')
\le \lambda_\varepsilon^n, 
\end{align*}
as required.
The theorem \ref{thm1} is proved. \hfill $\square$

\begin{remark}
In other words, the theorem \ref{thm1} claims that whatever exponential rate 
$$
\sup_x \|\mu^x_k - \mu_{inv}\|_{TV} \le C \lambda^k, \quad k\ge 1, 
$$
with some $C>0$ and $\lambda \in (0,1)$ 
holds true, this rate may be (nearly) attained by  the iterated Markov -- Dobrushin's method by choosing a sufficiently large number of iterations $m$. 

\end{remark}


\subsection{
Operators $V^{(m)}$ and 
$\hat V^{(m)}$ and their spectral radii}\label{sec:homo2} 
In this subsection, we present the bound, which is intermediate between Gantmacher and iterated Markov -- Dobrushin ones  in the finite state space case. The bounds themselves are valid in general state spaces. The author  believes that this intermediate estimate is the best available in the cases where Gantmacher's approach is not working\footnote{With a remark that, in principle, the step size $m$ may vary.}.


\subsection*{Coupling lemma}
Let us recall the folklore ``coupling'' lemma and the basic construction of coupling taken from \cite[\S 2.1]{VerVer}, not pretending to be an author of the lemma itself. The coupling approach was created by  W. Doeblin (V. Doblin) \cite{Doeblin}, who, in turn, referred to A.N. Kolmogorov \cite{Kolm}, even though, of course, there is no such term or notion in the latter seminal paper on Markov processes. 
In any case, the construction taken from \cite{Ver_Doeblin, VerVer} and its modification for $m$ steps will be used in what follows. Slightly different versions could be found in \cite{Thor, Kulik-Scheutzow}, and a large number of semi-Markov models in relation to coupling may be found in \cite{Silvestr}. We highlight that it is the idea of iterations based on the proposed coupling construction which leads to the (nearly) best possible rate of convergence; at least, it is certainly so for finite state spaces due to F.R. Gantmacher \cite[formula (13.96)]{Gant}. 

\begin{lemma}
\label{odvuh}
Let $X^{1}$ and $X^2$ be two random variables on their (without loss of generality different, which will be made independent after we take their direct product) probability spaces $(\Omega^1, {\cal F}^1, \mathbb P^1)$ and $(\Omega^2, {\cal F}^2, \mathbb P^2)$ and with densities $p^1$ and $p^2$ with respect to some reference measure $\Lambda$, correspondingly.  Then, if 
\begin{equation*}
\alpha := \int \left(p^1(x)\wedge p^2(x)\right) \Lambda(dx) > 0, 
\end{equation*}
then there exists one more probability space $(\Omega, {\cal F}, \mathsf P)$ and two random variables on it $\tilde X^1, \tilde X^2$ such that 
\begin{equation*}
{\cal L}(\tilde X^j) ={\cal L}(X^j), \; j=1,2, \qquad  \quad \mathsf  P(\tilde X^1 = \tilde X^2) = \alpha. 
\end{equation*}
\end{lemma}

\medskip

\noindent 
{\bf Construction and proof.} We will need now {\em four} new independent random variables, Bernoulli random variable  
$\zeta$ with $\mathsf P(\zeta=0) = \delta$ and $\eta^{1,2}$  and $\xi$ with the densities with respect to the measure $\Lambda$, respectively, 
\begin{eqnarray}\label{laws}
& \displaystyle p^{\eta^1}(x) := \frac{p^1 - p^1\wedge p^2}{\displaystyle\int (p^1 - p^1\wedge p^2)(y)\Lambda(dy)}(x), \quad
p^{\eta^2}(x) := \frac{p^2 - p^1\wedge p^2}{\displaystyle\int (p^2 - p^1\wedge p^2)(y)\Lambda(dy)}(x),
 \nonumber \\ \\ \nonumber 
& \displaystyle p^{\xi}(x) := \frac{ p^1\wedge p^2}{\displaystyle\int (p^1\wedge p^2)(y)\Lambda(dy)}(x),
\end{eqnarray}
where in the last expression it is assumed that the denominator is strictly  positive; the alternative case will be explained in the end of the proof; in the first two expressions it is also assumed that the denominator is strictly positive, and the alternative will be treated in the last step of the proof. 

We may assume that these four r.v.'s, $\zeta, \eta^{1}, \eta^2, \xi$ are all defined on their own probability spaces and eventually we consider the {\bf direct product} of these probability spaces denoted as $(\Omega, {\cal F}, \mathbb P)$. As a result, they are all defined on one unique probability space and are independent there. 
Now, on the same product of all  probability spaces just mentioned, let us define the random variables 
\begin{equation}\label{ety}
\tilde X^1:=  \eta^1 1(\zeta \not=0) + \xi 1(\zeta=0); \quad  
\quad \tilde X^2:=  \eta^2 1(\zeta \not=0) +\xi 1(\zeta=0).
\end{equation}
Then it holds that (see the detailed proof in  \cite{VerVer})
\begin{eqnarray*}
&\displaystyle \mathsf E g(\tilde X^i)  
= \mathsf E g(X^i), 
\quad i=1,2,
\end{eqnarray*}and also, 

\[
\mathsf P(\tilde X^1=\tilde X^2) = \alpha.  
\]
The detailed proof may be found in \cite{VerVer}. \hfill $\square$

\subsection*{Markov coupling, homogeneous case, $m$ steps.}
Let us now see how to generalize the Lemma \ref{odvuh} to a sequence of random variables and present our coupling
construction for homogeneous Markov chains based on \cite{VerVer}. Consider two versions $(X^1_n), (X^2_n)$ of the same Markov process with two initial distributions $\mu_0^1$ and  $\mu_0^2$ respectively (this does not exclude the case of non-random initial states). Denote 
\begin{equation*}
\alpha(0) :=
\int \left(\frac{\mu_0^1(dy)}{\mu_0^2(dy)}\wedge 1 \right)\mu_0^2(dy).
\end{equation*}
Notice that $0\le \alpha(0)\le1$ similarly to $\alpha(u,v)$ for all $u,v$.  
We assume that $X^1_0$ and $X^2_0$ have different distributions, so $\alpha(0)<1$. Otherwise, we would obviously have $X^1_n\stackrel{d}{=}X^2_n$ (equality in distribution) for all $n$, and then the coupling could have been made trivially, for example, by letting  $\widetilde X^1_n\!=\! \widetilde
X^2_n\!:\!=\!X^1_n$.

Let us introduce a new, vector-valued homogeneous {\bf Markov process} $\left(\eta^1_n,\eta^2_n,\xi_n,\zeta_n\right)$. Notations will be used like $x=(x^1,x^2,x^3,x^4)$, $y=(y^1,y^2,y^3,y^4)$, where all $x^i, y^i \in S$ for $1\le i\le 3$, and $x^4,y^4 = 0,1$. Hence, without any index, $x$ and $y$ are both four component vectors, while their components are denoted by indices like $x^1$, et al.  
The values $\left(\eta^1_0,\eta^2_0,\xi_0,\zeta_0\right)$ are chosen directly on the basis of the Lemma \ref{odvuh} as $\left(\eta^1,\eta^2,\xi,\zeta\right)$, according to the distributions in (\ref{laws}).
In particular, if $\delta_0=0$ then we can set
\begin{equation*}
\eta^1_0:=X^1_0,\; \eta^2_0:=X^2_0,\; \xi_0:=X^1_0,\; \zeta_0:=1.
\end{equation*}
(The value for $\xi_0$ is not important in this case.)
If $\delta_0=1$ then we can set
\begin{equation*}
\eta^1_0:=X^1_0,\; \eta^2_0:=X^1_0,\; \xi_0:=X^1_0,\; \zeta_0:=0.
\end{equation*}

Now, by induction, assuming that the random variables $\left(\eta^1_{nm},\eta^2_{nm},\xi_{nm},\zeta_{nm}\right)$ have been determined for some $n$, let us show how to construct them for $(n+1)m$. For this aim, we define the transition probability density $\phi$ with respect to the 
measure $\Lambda^{(m)}_{x^1, x^2} \times \Lambda^{(m)}_{x^1, x^2} \times \Lambda^{(m)}_{x^1, x^2}\times (\delta_0 + \delta_1)$ for this vector-valued process as follows,
\begin{equation}\label{process_eta}
\phi^{(m)}(x,y):=\phi^{(m)}_1(x,y^1)\phi^{(m)}_2(x,y^2)\phi^{(m)}_3(x,y^3) \phi^{(m)}_4(x,y^4),
\end{equation}
where, recall, $x\!=\!(x^1,x^2,x^3,x^4)$, $y\!=\!(y^1,y^2,y^3,y^4)$, and if
$0\!<\!\delta^{(m)}(x^1,x^2)\!<\!1$, then
\begin{align}
&\displaystyle \phi^{(m)}_1(x,u):=\frac{p^{(m)}(x^1,u)-p^{(m)}(x^1,
u)\wedge p^{(m)}(x^2,u)}{1-\delta^{(m)}(x^1,x^2)},\label{phi_11} 
 \\\nonumber\\ 
&\phi^{(m)}_2(x,u):=\frac{p^{(m)}(x^2,u)-p^{(m)}(x^1,u)\wedge p^{(m)}(x^2,u)}{1-\delta^{(m)}(x^1,x^2)},\label{phi_12}
 \\\nonumber\\
&\displaystyle \phi^{(m)}_3(x,u):=1(x^4=1)\frac{p^{(m)}(x^1,u)\wedge
 p^{(m)}(x^2,u)}{\delta^{(m)}(x^1,x^2)}+1(x^4=0)p^{(m)}(x^3,u), \label{phi_3}
 \\\nonumber\\ \nonumber
&\displaystyle \phi^{(m)}_4(x,u):=1(x^4=1)\left(\delta_1(u)(1-\delta^{(m)}(x^1,x^2))+ 
\delta_0(u)\delta^{(m)}(x^1,x^2)\right) 
 \\ \nonumber
 \\ 
&\displaystyle 
\hspace{2cm} +1(x^4=0)\delta_0(u). \label{phi_4}
\end{align}
Here $\delta_i(u)$ is the Kronecker symbol, $\delta_i(u) = 1(u=i)$, or, in other words, the delta measure concentrated at state $i$. The case $x^4=0$ signifies coupling which has already been realised at the previous step, and $u=0$ means successful coupling at the transition.  
Note that $\phi_1$ and $\phi_2$ do not depend on the variable $x^3$; we will denote it by the notation $\phi_i((x^1,x^2,*,x^4),u)$ ($i=1,2$) where $*$ stands for any possible value of $x^3$. Also even if it is written $\phi_3((x^1,x^2,x^3,1),u)$, yet, this value does not depend on $x^3$ either.

In the degenerate cases, if $\delta^{(m)}(x^1,x^2)=0$ (coupling at the transition is impossible), then instead of (\ref{phi_3}) we set,  e.g., 
$$
\phi^{(m)}_3(x,u):=1(x^4=1)p(x^3,u) + 1(x^4=0)p^{(m)}(x^3,u) = p^{(m)}(x^3,u),
$$ 
and if $\delta^{(m)}(x^1,x^2)=1$, then instead of (\ref{phi_12}) we may set 
$$
\phi^{(m)}_1(x,u)=\phi^{(m)}_2(x,u):= p^{(m)}(x^1,u). 
$$ 
However, notice that in the case of $\delta^{(m)}(x^1,x^2)=1$ we shall not assume that the ``next'' values $\eta^1_{(n+1)m}$ and  $\eta^2_{(n+1)m}$ for any $n\ge 1$ are not equal, simply  because this event may only occur with probability zero. This is not a problem, but should be remembered of in the Markov coupling construction. 
The formula (\ref{phi_4}) which defines \(\phi^{(m)}_4(x,u)\) can be accepted in all cases.

\medskip

Finally, let 
\begin{align}\label{xin0}
\widetilde X^1_{nm}\!:=\!\eta^1_{nm} 1(\zeta_{nm}\!=\!1)\!+\!\xi_{nm} 1(\zeta_{nm}\!=\!0), \; 
\widetilde X^2_{nm}\!:=\!\eta^2_{nm} 1(\zeta_{nm}\!=\!1)\!+\!\xi_{nm} 1(\zeta_{nm}\!=\!0).
\end{align}

It may be verified (see \cite{VerVer}) that the transition  densities for the components $\tilde X^1_{(n+1)m}$ and $\tilde X^2_{(n+1)m}$, respectively, given $\tilde X^1_{nm}$ and $\tilde X^2_{nm}$ depend functionally not on both $\tilde X^1_{nm}$ and $\tilde X^2_{nm}$,  but  the first one depends (functionally) only on $\tilde X^1_{nm}$, and the second one, respectively, on $\tilde X^2_{nm}$. 
It may be checked (see again \cite{VerVer}) that, due to the construction above, we have, in particular, the following densities of the conditional distributions of $(\tilde X^1_{(n+1)m}, \tilde X^2_{(n+1)m})$ given $(\tilde X^1_{nm}, \tilde X^2_{nm})$:  
\begin{align*}
&\frac{\mathsf P(\tilde X^1_{(n+1)m} \in dx^1 |\tilde X^1_{nm}, \tilde X^2_{nm})}{\Lambda_{\tilde X^1_{nm}, \tilde X^2_{nm}}(dx^1)} 
\!=\! p^{(m)}(\tilde X^1_{nm},x^1), 
 \\\\
&\frac{\mathsf P(\tilde X^2_{(n+1)m} \in dx^2 |\tilde X^1_{nm}, \tilde X^2_{nm})}{\Lambda_{\tilde X^1_{nm}, \tilde X^2_{nm}}(dx^2)} 
= p^{(m)}(\tilde X^2_{nm},x^2),
\end{align*}
in all cases.

Due to all of these, each of the components $(\tilde X^1_{nm}, \, n\ge 0)$ and $(\tilde X^2_{nm}, \, n\ge 0)$ are Markov processes with the same generator as $(X^1_{nm}, \, n\ge 0)$ and $(X^2_{nm}, \, n \ge 0)$. Recall that here $m$ is fixed, while $n$ runs from $0$ to infinity, and the common generator for both processes is $L^{(m)}h(x):= \mathsf E_x h(\tilde X^1_m) - h(x)$, although, the notation $L^{(m)}$ will not be used in what follows.
Moreover, the following lemma holds true. 
\begin{lemma}\label{lemma:2} 
Let the random variables $\widetilde X^1_{nm}$ and $\widetilde X^2_{nm}$, for $n\in\mathbb{Z}_+$ be defined   by the formulae:
\begin{equation}\label{xin}
\tilde X^1_n:=\eta^1_n 1(\zeta_n=1)+\xi_n 1(\zeta_n=0), \quad 
\tilde X^2_n:=\eta^2_n 1(\zeta_n=1)+\xi_n 1(\zeta_n=0).
\end{equation}
Then 
\begin{equation}\label{l71}
\widetilde X^1_{nm}\stackrel{d}{=}X^1_{nm}, \;\;\widetilde
 X^2_{nm}\stackrel{d}{=}X^2_{nm}, \quad \mbox{for all $n\ge 0$,}
\end{equation}
which implies that the process $\widetilde X^1_{nm}$ is equivalent to $X^1_{nm}$, and the process $\widetilde X^2_{nm}$ is equivalent to $X^2_{nm}$ in distribution in the space of trajectories; in particular, each of them is a Markov process with the same generator as $X^1_{nm}$.
Moreover, the couple $\tilde X_{nm}:=\left(\widetilde X^1_{nm}, \widetilde X^2_{nm}\right)$, $n\ge 0$, is also  a  homogeneous Markov process, and 
\begin{equation*}\label{mpequi}
\left(\widetilde X^1_{nm}\right)_{n\ge 0}\stackrel{d}{=}\left(X^1_{nm}\right)_{n\ge 0}; \quad \quad 
\left(\widetilde X^2_{nm}\right)_{n\ge 0}\stackrel{d}{=}\left(X^2_{nm}\right)_{n\ge 0}.
\end{equation*}
Also, 
\begin{equation}\label{l72}
\widetilde X^1_n=\widetilde X^2_n, \quad \forall \; n\ge
n_0(\omega):=\inf\{k\ge0: \zeta_{km}=0\}, 
\end{equation}
and
\begin{equation}\label{estimate}
\frac12\, \|\mu^{x^1}_{nm} - \mu^{\mu_0}_{nm}\|_{TV}\le 
\P_{x^1,\mu_0}(\widetilde X^1_{nm}\neq \widetilde
 X^2_{nm})\le \E_{x^1,\mu_0}\prod_{i=0}^{n-1}
 (1-\delta^{(m)}(\eta^1_i,\eta^2_i)).
\end{equation}
\end{lemma}
The detailed proof in the case $m=1$ may be found in \cite{VerVer}; for general $m$ the changes are negligible. \hfill $\square$

\medskip

The non-homogeneous  version of this lemma is offered  in what follows. 

\begin{remark}
The reader may ask, why this construction only provides the random variables at times $(mn, \, n\ge 1)$ where $m$ is fixed, and not for all $n\ge 1$ values. The matter is that as far as convergence, mixing, or coupling is estimated for some subsequence $(X_{n_i}, \,n_i\to\infty)$, in the case of  Markov chains it is immediately inherited by the whole process $(X_{n}, n\ge 1)$, perhaps with some logarithmically negligible multiplier (i.e., uniformly bounded).

\end{remark}



\subsection*{Operators $V^{(m)}$ and $\hat V^{(m)}$ and their spectral radii}\label{sec:Vr}

Now, let 

\begin{equation}\label{V}
V^{(m)}h(x) := (1-\delta^{(m)}(x_1,x_2)) \mathsf E_{x_1,x_2}h(\tilde X_m), 
\end{equation}
for any bounded, Borel measurable  function $h$ on the space $S^2 := S\times S$, and  with $x=(x_1, x_2)\in S^2$; recall that 
$\tilde X_{nm} = (\tilde X_{nm}^1, \tilde X^2_{nm})$. Notice that on the diagonal $(x=(x_1,x_2) : \, x_1=x_2)$ we have 
$$
V^{(m)}h(x) = (1-\delta^{(m)}(x_1,x_1))\mathsf E_{x_1,x_2}h(\tilde X_m)\vert_{x_2=x_1} = 0. 
$$
Hence, it makes sense to reduce the operator itself on functions defined on 
$\hat S^2 := S^2 \setminus \text{diag}(S^2)$, 
that is, to define for $x=(x_1,x_2)\in \hat S^2$ and for functions $\hat h: \hat S^2 \to \mathbb R$, 
\begin{equation}\label{Vhat}
\hat V^{(m)}\hat h(x):= (1-\delta^{(m)}(x_1,x_2))\mathsf E_{x^1,x^2}\hat h(\tilde X_m).
\end{equation}

The estimate (\ref{estimate}) can be rewritten via the operator $V^{(m)}$, or, equivalently,  via $\hat V^{(m)}$ as follows: 
\begin{align}\label{l2}
\!\!\P_{x_1,\mu}(\widetilde X^1_{nm}\neq \widetilde X^2_{nm})\!\le\! 
\int \E_{x_1,x_2} (V^{(m)})^n {\bf 1}(x_1,x_2)1(x_1\!\not =\!x_2)\mu(dx_2) 
 \nonumber\\  \\\nonumber
=\int \E_{x_1, x_2}1(\tilde X_0^1\!\not = \!\tilde X_0^2) (\hat V^{(m)})^n {\bf 1}(x_1,x_2)1(x_1\!\not =\!x_2)\mu(dx_2). 
\end{align}
Here in the first line of (\ref{l2}), $\bf 1$ is the function on $S^2$ identically equal to one, and $(V^{(m)})^n {\bf 1}$ is the result of the operator $(V^{(m)})^n$, which acts on this function, like $(V^{(m)})^n h$ with function $h(x) = {\bf 1}(x)$. In turn, $(V^{(m)})^n {\bf 1}(x_1,x_2)$ is the value of the function $(V^{(m)})^n {\bf 1}$ at $x=(x_1,x_2)$. In the second line of (\ref{l2}), $\bf 1$ is the function identically equal to one on $\hat S^2$. So, {\bf we obtain the following}

\begin{theorem}\label{cor} 
Let the homogeneous Markov chain $(X_n)$ be irreducible and aperiodic. Then
for any fixed $m$ the estimate holds, 
\begin{equation}\label{estimate2}
\frac12\, \|\mu^{x^1}_{nm} - \mu^{x^2}_{nm}\|_{TV}\le \|V^{(m)}\|^n, \quad \forall \, n\ge 1.
\end{equation}
\end{theorem}

\noindent
Now, the well-known inequality (see, for example, \cite[\S 8]{KLS}) reads, 
\begin{equation*}\label{rkappa}
r(V^{(m)}) \le \|V^{(m)}\|. 
\end{equation*}
Further, notice that since the operator $V^{(m)}$ is positive (that is, it transforms any non-negative function again to non-negative), we have,  
\begin{equation}\label{Vmam}
\|V^{(m)}\| = \sup\limits_{x\in S^2} V{\bf 1}(x) = \sup_{x\in S^2} (1-\alpha^{(m)}(x)) = 1-\alpha^{(m)}.
\end{equation} 

The same holds true for the operator $\hat V^{(m)}$ (here the function ${\bf 1}(x)\equiv 1$ is defined on $x\in \hat S^2$): 
\begin{equation*}\label{rkappa}
r(\hat V^{(m)}) \le \|\hat V^{(m)}\| = \sup\limits_{x\in \hat S^2} \hat V{\bf 1}(x) = \sup_{x\in \hat S^2} (1-\alpha^{(m)}(x)) = 1-\alpha^{(m)}. 
\end{equation*}

It will be easier to argue with the operator $V^{(m)}$ in the sequel; so, we will continue with  this operator. However, from the computational point of view -- that is, to compute the spectral radius -- the operator $\hat V^{(m)}$ may be preferred because of some reduction in dimension. Such a reduction is possible because for any function $h(x), x=(x_1,x_2)$, on the diagonal $x^1 = x^2$ we have $V h(x)=0$, and because 
$\hat V^{(m)}$ is the projection of $V^{(m)}$ on $\hat S^2$. We do not pursue further this issue here.

Because the operator $V^{(m)}$ is positive, we have, 
\[
\lim_{n \to\infty} \frac1n \, \ln (V^{(m)})^n {\bf 1}(x) \le \lim_{n \to\infty} \frac1n \, \ln \|(V^{(m)})^n\|. 
\]
The well-known Gelfand formula for the spectral radius (see, for example, \cite[formula (8.4)]{KLS}) claims that  \(r(V^{(m)})= \lim_{n\to\infty}\|(V^{(m)})^n\|^{1/n}\), that is, 
\[
\lim_{n \to\infty} \frac1n \, \ln \|(V^{(m)})^n\| = \ln r(V^{(m)}).
\]
Also, recall that always
\[
r(V^{(m)}) \le \|V^{(m)}\|.
\]
As a consequence, and due to the theorem \ref{thm1}, we obtain the following result. 
\begin{theorem}\label{thm2} 
Let the homogeneous Markov chain $(X_n)$ be irreducible and aperiodic. Then the following assertions hold true.\\
1. In all cases, for any $x^1 \in S$ and any $m\ge 1$,
\begin{align}\label{newrate12}
&\limsup\limits_{n\to\infty} \frac1{n} \ln \| \mu^{x_1}_n - \mu_{inv}\|_{TV} 
\!\le\! \limsup_{n\to\infty} \frac1{nm}  \ln\! \int \!\!2(V^{(m)})^n {\bf 1}(x_1, x_2)\mu(dx_2)
 \nonumber \\\\ \nonumber
&\!\le\! \frac1{m} \ln r(V^{(m)}))^{1/m}. 
\end{align}

\medskip

\noindent
2. If there exists $\lambda>0$ such that for any $\varepsilon>0$ there exist $C$ such that 
\[
\sup_{x_1,x_2} \|\mu^{x_1}_m - \mu^{x_2}_m\|_{TV}\le C(\lambda+\varepsilon)^m 
\]
for some $m\ge 1$ 
(see (\ref{thm1-e1})), then 
\begin{equation}
r(V^{(m)}) \le \lambda+\varepsilon.
\end{equation}
\end{theorem}
Recall that often the role of such a $\lambda$ is played by $|\lambda_2({\cal P})|$, see \cite[formula (13.96)]{Gant}.

\begin{corollary}\label{corr1} 
Let the assumptions of the theorem \ref{thm2} be satisfied.
Then, under the condition
\begin{equation}\label{r12}
r(V^{(m)})<1, 
\end{equation}
the rate of convergence in
\[
\|\mu^x_n - \mu_{inv}\|_{TV} \to 0, \quad n\to\infty
\]
is exponential: 
for any $\varepsilon>0$ and $n$ large enough ($n\ge N(x,\varepsilon)$), 
\begin{equation}\label{newrate2}
\|\mu^x_n - \mu_{inv}\|_{TV}
\le (r(V^{(m)})+ \varepsilon)^{[n/m]}.  
\end{equation}
\end{corollary}

\begin{corollary}\label{corr2} 
Let the assumptions of the theorem \ref{thm2} be satisfied. 
Then, in the case of finite $S$ the equality holds true,
$$
|\lambda_2({\cal P})| = \lim_{m\to\infty} \big(r(V^{(m)})\big)^{1/m}.
$$
\end{corollary}

\begin{remark}
The bound (\ref{newrate2}) is better in a non-strict sense  than (\ref{estimate2}). They both ``nearly'' attain the rate provided by the second eigenvalue $\lambda_2({\cal P})$, if $m$ is large enough. 
Let us recall some observations from \cite[\S 4]{VerVer}. They all concern comparisons of three methods for particular transition matrices: $\lambda_2({\cal P})$, simple MD, and simple $r(V)$. For very simple examples with small cardinality of $|S|$, all three methods may provide the same asymptotics. For larger values of $|S|$, ``usually'' the method by $r(V)$ provides strictly better estimate than a ``simple'' MD-1 method, and the $\lambda_2({\cal P})$ is usually considerably better than both $r(V)$ and MD-1. Yet, the authors did not compare there $\lambda_2({\cal P})$ with $r(V^{(m)})$ and MD-$m$ with $m>1$.
\end{remark}

\begin{remark}
In general, 
$$
r(V^{(m)})^{1/m}\, \neq \, r(V), 
$$
because, generally speaking, $V^{(m)} \neq V^m$, as may be seen from the examples in \cite{NV1}. So, the remark 13 in \cite{VerVer}, which claims (without proof) the opposite 
is wrong. This does not affect any of the results in \cite{VerVer}, except for the remark itself. 

\end{remark}

\section{Non-homogeneous case}
The non-homogeneous version may be  developed  in a similar way. The major difference is, of course, that in this case there is neither second eigenvalue approach, nor spectral radius of $V$, or $V^{(m)}$, since in the non-homogeneous case all operators may depend on time, and so may be different and non-commutative. The remaining tools are the non-homogeneous Markov -- Dobrushin method, its iterations and their improvements based on the norms of the products $\|\prod_{}^{} V_k\|$ and $\|\prod_{k}^{} V^{(m)}_{mk}\|$. 
In \cite{VerVer} the case $m>1$ was just outlined. Hence, more details will be now presented. 
Let \begin{equation}\label{MDt}
 \alpha_t^{(m)} := \inf_{x_1,x_2} \int \left(\frac{\P_{t,x_2}(t+m,dy)}{\P_{t,x_1}(t+m,dy)}\wedge 1 \right)\P_{t,x_1}(t+m,dy), 
\end{equation}
and 
$$
\alpha_t^{(m)}(x_1,x_2) :=  \int \left(\frac{\P_{t,x_2}(t+m,dy)}{\P_{t,x_1}(t+m,dy)}\wedge 1 \right)\P_{t,x_1}(t+m,dy). 
$$
Further, for any $m\ge 1$ and for any two or three fixed states $x_1, x_2, x_3$ and for any $t\ge 0$ denote 
\begin{align*}
&\Lambda^{(m)}_{t,t+m,x_1,x_2}(dz) := (\P_{t,x_1}(t+m,dz) + \P_{t,x_2}(t+m,dz))/2,
 \\
&\Lambda^{(m)}_{t,t+m,x_1,x_2, x_3}(dz) := (\P_{t,x_1}(t+m,dz) + \P_{t,x_2}(t+m,dz) + \P_{t,x_3}(t+m,dz))/3.
\end{align*} 
\noindent
Notice that $\Lambda_{t,x_1,x_2}(dz) = \Lambda_{t,x_1,x_2}(dz)$, and similarly $\Lambda_{t,x_1,x_2,x_3}$ as well as $\Lambda^{(m)}_{t,x_1,x_2,x_3}$ do not depend on the permutation of the variables $(x_1,x_2,x_3)$.
\begin{lemma}\label{newMDt} 
The following representation for the coefficient $\alpha_t^{(m)}$ in (\ref{MDt}) holds true, 
\begin{equation}\label{betterMD}
 \alpha^{(m)}_t = \inf_{x_1,x_2} \int \left(\frac{\P_{t,x_2}(t+m,dy)}{\Lambda^{(m)}_{t,x_1,x_2}(dy)}\wedge \frac{\P_{t,x_1}(t+m,dy)}{\Lambda^{(m)}_{t,x_1,x_2}(dy)} \right)\Lambda^{(m)}_{t,x_1,x_2} (dy).
\end{equation}
In particular -- since $\Lambda_{t,x_1,x_2}(dz) = \Lambda_{t,x_2,x_1}(dz)$ -- 
for any $x_1,x_2\in S$ we have,
\[
\alpha^{(m)}_t(x_1,x_2) = \alpha^{(m)}_t(x_2,x_1).
\]
\end{lemma}

\noindent
{\bf Proof} is straightforward. \hfill $\square$

\medskip

Similarly to the homogeneous case, we obtain the following. 
\begin{proposition}\label{pro3} 
The uniform bound holds true for any $n$ 
and for any $m\ge 1$,
\begin{equation}\label{pro3-e1}
\sup_{x_1,x_2}\sup_{A\subset S} |\mu^{x_1}_n(A) - \mu^{x_2}_n(A)| \le   \Big(\prod_{t=0}^{[(n-1)/m]} (1-\alpha_{tm}^{(m)})\Big) (1-\alpha_{[n/m]m}^{(n - [n/m]m)}), 
\end{equation}
and 
\begin{equation}\label{pro3-e2}
\sup_{x_1,x_2}\|\mu^{x_1}_n - \mu^{x_2}_n\|_{TV} \le 2  \Big(\prod_{t=0}^{[(n-1)/m]} (1-\alpha_{tm}^{(m)})\Big) (1-\alpha_{[n/m]m}^{(n - [n/m]m)}). 
\end{equation}
{\bf Proof}  follows easily from the construction and Markov property. \hfill $\square$

\end{proposition}
\begin{remark}
These two estimates are provisional and may be regarded as some extensions of Kolmogorov's inequality \cite[formula (29)]{Kolm}; the latter corresponds to the case $m=1$. More than that, in fact, a bound  similar\footnote{Actually, it even allows non-equal differences between times $0<m_1<m_2<\ldots$ in place of $0<m<2m<\ldots$, which is, clearly, also possible in a bound analogous to  (\ref{pro3-e1}), as well as in the theorem \ref{thm3} and in  its corollary in what follows.} to our inequality (\ref{pro3-e1}) can be found -- without a number -- in the proof of the theorem 3 in \cite[\S 4]{Kolm}. So far, the only moderate development in comparison to these bounds in \cite[\S 4]{Kolm} is that we know that in the homogeneous case, the larger is $m$, the closer the bound is to the best possible one, which is due to the Gantmacher result mentioned earlier (see a simplified version in the proposition \ref{pro2}). Still, the application of coupling technique provides, in principle, more precise bounds than (\ref{pro3-e1}) and (\ref{pro3-e2}) for any finite value of $m$, see the theorem \ref{thm3} and the corollary \ref{cor3} in what follows.
\end{remark}


\subsection*{Markov coupling, non-homogeneous case}
As earlier, denote
\begin{equation*}
\alpha(0) :=
\int \left(\frac{\mu_0^1(dy)}{\mu_0^2(dy)}\wedge 1 \right)\mu_0^2(dy).
\end{equation*}
The random variables $\left(\eta^1_0,\eta^2_0,\xi_0,\zeta_0\right)$ are chosen using the lemma \ref{odvuh} as $\left(\eta^1,\eta^2,\xi,\zeta\right)$, according to the distributions in (\ref{laws}), exactly as in the homogeneous case. Recall that we fix $m\ge 1$.

Further, by induction, assuming that the random variables $\left(\eta^1_{nm},\eta^2_{nm},\xi_{nm},\zeta_{nm}\right)$ have been determined for some $n$, let us show how to construct them for $m(n+1)$. For this aim, let us define the transition probability density $\phi^{(m)}_n$ with respect to the measure 
$\Lambda^{(m)}_{(n+1)m,x^1, x^2} \times \Lambda^{(m)}_{(n+1)m,x^1, x^2} \times \Lambda^{(m)}_{(n+1)m,x^1, x^2}\times (\delta_0 + \delta_1)/2$ for this vector-valued process as follows (in what follows $t=nm$, where $m\ge 1$ is fixed and $n = 1,2,\ldots$),
\begin{equation}\label{process_eta2}
\phi^{(m)}_t(x,y):=\phi^{(m)}_{1,t}(x,y^1)\phi^{(m)}_{2,t}(x,y^2)\phi^{(m)}_{3,t}(x,y^3) \phi^{(m)}_{4,t}(x,y^4),
\end{equation}
where $x=(x^1,x^2,x^3,x^4)$, $y=(y^1,y^2,y^3,y^4)$, and if
 $0<\alpha_t(x^1,x^2)<1$, then
\begin{align}
&\displaystyle \phi^{(m)}_{t,1}(x,u):=\frac{p^{(m)}_t(x^1,u)-p^{(m)}_t(x^1,
u)\wedge p^{(m)}_t(x^2,u)}{1-\alpha^{(m)}_t(x^1,x^2)},  
\label{phi_1211x1}
 \\\nonumber\\
&\phi^{(m)}_{t,2}(x,u):=\frac{p^{(m)}_t(x^2,u)-p^{(m)}_t(x^1,u)\wedge p^{(m)}_t(x^2,u)}{1-\alpha^{(m)}_t(x^1,x^2)},\label{phi_1211x2}
 \\\nonumber\\
&\displaystyle \phi^{(m)}_{t,3}(x,u):=1(x^4=1)\frac{p^{(m)}_t(x^1,u)\wedge
 p^{(m)}_t(x^2,u)}{\alpha^{(m)}_t(x^1,x^2)}+1(x^4=0)p^{(m)}_t(x^3,u),\label{phi_311}
 \\\nonumber\\ 
&\displaystyle 
\phi^{(m)}_{t,4}(x,u):=1(x^4=1)\big(\delta_1(u)(1-\alpha^{(m)}_t(x^1,x^2))+ 
\delta_0(u)\alpha^{(m)}_t(x^1,x^2)\big)
  \nonumber \\ \nonumber\\
& 
\hspace{2cm}+1(x^4=0)\delta_0(u)\label{phi_411}. 
\end{align}
The case $x^4=0$ signifies coupling which has already been realised at the previous step, and $u=0$ means successful coupling at the transition.  
More details and comments may be found in \cite{VerVer}, mostly about the case $m=1$.

Let the random variables $\widetilde X^1_{nm}$ and $\widetilde X^2_{nm}$, for $n\in\mathbb{Z}_+$ be defined   by the following formulae:
\begin{align}\label{xin3}
\widetilde X^1_{nm}:=\eta^1_{nm} 1(\zeta_{nm}=1)+\xi_{nm} 1(\zeta_{nm}=0), 
 \nonumber \\\\ \nonumber
\widetilde X^2_{nm}:=\eta^2_{nm} 1(\zeta_{nm}=1)+\xi_{nm} 1(\zeta_{nm}=0).
\end{align}
Then, similarly to the homogeneous case, it may be shown that
\begin{align*}
\frac{\mathsf P(\tilde X^1_{(n+1)m} \in dx^1 |\tilde X^1_{nm}, \tilde X^2_{nm})}{\Lambda_{(n+1)m,\tilde X^1_{nm}, \tilde X^2_{nm}}(dx^1)} 
= p_{nm}(\tilde X^1_{nm},x^1),
\end{align*}
and
\begin{align*}
\frac{\mathsf P(\tilde X^2_{(n+1)m} \in dx^2 |\tilde X^1_{nm}, \tilde X^2_{nm})}{\Lambda_{(n+1)m,\tilde X^1_{nm}, \tilde X^2_{nm}}(dx^2)} 
= p_{nm}(\tilde X^2_{nm},x^2),
\end{align*}
in all cases. This signifies that each component $\tilde X^i_{nm}$ is a Markov process equivalent to $X^i_{nm}$, $i=1,2$ (cf.  to \cite[lemma 18]{VerVer}). Similarly to the case $m=1$ as in \cite{VerVer}, it may be shown that the pair $(X^1_{nm}, X^2_{nm})$ is also a Markov process and that the following lemma holds true. 

\begin{lemma}\label{lem4} 
For any $m\ge 1$, 
\begin{equation}\label{l71}
\widetilde X^1_{nm}\stackrel{d}{=}X^1_{nm}, \;\;\widetilde
 X^2_{nm}\stackrel{d}{=}X^2_{nm}, \quad \mbox{for all $n\ge 0$,}
\end{equation}
which implies that the process $\widetilde X^1_{nm}$ is equivalent to $X^1_{nm}$, and the process $\widetilde X^2_{nm}$ is equivalent to $X^2_{nm}$ in distribution in the space of trajectories; in particular, each of them is a Markov process with the same generator as $X^1_{nm}$.
Moreover, the couple $\tilde X_{nm}:=\left(\widetilde X^1_{nm}, \widetilde X^2_{nm}\right)$, $n\ge 0$, is also  a  (non-homogeneous) Markov process with the transition density
with respect to $\Lambda_{(n+1)m,\tilde X^1_{nm}, \tilde X^2_{nm}}$,
\begin{align*}
&\frac{\mathsf P(\tilde X^1_{(n+1)m} \in dx^1, \tilde X^2_{(n+1)m} \in dx^2 |\tilde X^1_{nm}, \tilde X^2_{nm})}{\Lambda_{(n+1)m,\tilde X^1_{nm}, \tilde X^2_{nm}}(dx^1)}
 \\\\
&\hspace{4mm}= 1(\tilde X^1_{nm} \not =\tilde X^2_{nm}) p_{nm}(\tilde X^1_{nm}, x^1) p_{nm}(\tilde X^2_{nm},x^2) 
 \\\\
&\hspace{4mm}+1(\tilde X^1_{nm} =\tilde X^2_{nm}) p_{nm}(\tilde X^1_{nm}, x^1)\delta(x_1-x_2).
\end{align*}
Also,
\begin{equation*}\label{mpequi}
\left(\widetilde X^1_{nm}\right)_{n\ge 0}\stackrel{d}{=}\left(X^1_{nm}\right)_{n\ge 0}; \quad \quad 
\left(\widetilde X^2_{nm}\right)_{n\ge 0}\stackrel{d}{=}\left(X^2_{nm}\right)_{n\ge 0}.
\end{equation*}
Moreover, 
\begin{equation}\label{l723}
\widetilde X^1_{nm}=\widetilde X^2_{nm}, \quad \forall \; n\ge
n_0(\omega):=\inf\{k\ge0: \zeta_k=0\}, 
\end{equation}
and
\begin{equation}\label{lem4-e3}
\P_{\mu^1,\mu^2}(\widetilde X^1_{nm}\neq \widetilde
 X^2_{nm})\le
\E_{\mu^1,\mu^2}\Big(\prod_{i=0}^{n-1}
(1-\alpha^{(m)}_{im}(\eta^1_{im},\eta^2_{im}))\Big). 
\end{equation}
\end{lemma}
{\bf Proof.} The claim was proved in \cite[lemma 18]{VerVer} for the case of $m=1$. The proof for a general $m$ is practically identical. The bound (\ref{lem4-e3}) may be regarded as a further extended version of Kolmogorov's estimate mentioned above. \hfill $\square$

\subsection*{Operators $V^{(m)}$ and $\hat V^{(m)}$ in the non-homogeneous case}\label{sec:Vr}

Let us introduce the operators $V^{(m)}_t$ acting on bounded, Borel measurable functions $h$ on the space $S^2 := S\times S$ as follows: for $x=(x_1, x_2)\in S^2$,
\begin{equation}\label{V3}
V^{(m)}_t h(x) := (1-\alpha^{(m)}_t(x_1,x_2)) \mathsf E_{t,x_1,x_2}h(\tilde X_{t+m}), 
\end{equation}
where 
$\tilde X_n = (\tilde X_n^1, \tilde X^2_n)$. Note that on the diagonal $x=(x_1,x_2) : x^1=x^2$ we have 
$$
V^{(m)}_th(x) = (1-\alpha^{(m)}_t(x_1,x_1))\mathsf E_{t,x_1,x_2}h(\tilde X_{t+m}) \vert_{x_2 = x_1}= 0, 
$$
since $\alpha^{(m)}_t(x_1,x_1) = 1$ for any $x_1\in S$. Hence, similarly to the homogeneous case, it makes sense to either consider the functions $h$ on $S^2$ vanishing on the diagonal $\text{diag}(S^2)= (x =  (x_1,x_1)\in S^2)$, or, equivalently, to reduce the operator itself to functions defined on 
$$
\hat S^2 := S^2 \setminus \text{diag}(S^2), 
$$
that is, to define for $x=(x^1,x^2)\in \hat S^2$ and for functions $\hat h: \hat S^2 \to \mathbb R$, 
$$
\hat V^{(m)}_t\hat h(x):= (1-\alpha^{(m)}_t(x^1,x^2))\mathsf E_{t,x^1,x^2}\hat h(\tilde X_{t+m}). 
$$
Notice that 
$$
\| \hat V^{(m)}_k\| = \|V^{(m)}_k\|.
$$

\begin{lemma}\label{lem5} 
The bound (\ref{pro3-e1}) may be equivalently rewritten as
\begin{equation}\label{lem5-e1}
\sup_{x_1,x_2}\sup_{A\subset S} |\mu^{x_1}_n(A) - \mu^{x_2}_n(A)| \le   \Big(\prod_{k=0}^{[(n-1)/m]} \| V^{(m)}_{km}\|\Big) \| V^{(n - m[n/m])}_{m[n/m]}\|. 
\end{equation}
\end{lemma}
{\bf Proof} follows from the equality (\ref{Vmam}) applied to each of the operators in the right hand side of (\ref{lem5-e1}). \hfill $\square$

\medskip

Further, the estimate (\ref{lem4-e3}) can be rewritten via the operators $V_t$ (or, equivalently,  via $\hat V_t$) as follows: 
\begin{align*}
\!\!\P_{\mu^1,\mu^2}(\widetilde X^1_{nm}\neq \widetilde X^2_{nm})\!\le\! 
\iint 
\prod_{i=0}^{n-1}  V_i^{(m)}{\bf 1}(x_1,x_2)1(x_1\!\not =\!x_2)\,\mu^1(dx_1)\mu^2(dx_2).
\end{align*}

\begin{remark}
In some cases it might be more convenient in examples to use the ``preliminary'' bound
\begin{align*}
\P(\widetilde X^1_{nm}\neq \widetilde X^2_{nm}) \le  
\P (\zeta_{nm} = 1) \equiv \E \Big(\prod_{i=0}^{n-1}1(\tilde X^1_{im} \neq \tilde X^2_{im})\Big). 
\end{align*}

\end{remark}

\begin{theorem}\label{thm3} 
For any Markov chain, with any $m\ge 1$,
\begin{align}\label{thm3-e1}
&\limsup\limits_{n\to\infty} \frac1n \ln \| \P_{\mu^1}(n,\cdot) - \P_{\mu^2}(n,\cdot)\|_{TV}
  \nonumber \\\\ \nonumber
&\!\le\! \limsup\limits_{n\to\infty} \frac1{nm} \ln \iint 
\prod_{i=0}^{n-1} V^{(m)}_i{\bf 1}(x^1,x^2)1(x^1\!\not =\!x^2)\mu^1(dx_1) \mu^2(dx_2).
\end{align}
\end{theorem}

{\bf Proof} follows straightforwardly from the bound (\ref{lem4-e3}) and from the definition (\ref{V3}). \hfill $\square$

\begin{corollary}\label{cor3}
In all cases, 
\begin{align}\label{cor3-e1}
&\limsup\limits_{n\to\infty} \frac1{n} \ln \| \P_{\mu^1}(n,\cdot) - \P_{\mu^2}(n,\cdot)\|_{TV}
 \nonumber \\\\ \nonumber 
&\!\le\! 
\frac1{nm} \ln \iint 
\prod_{i=0}^{n-1} V^{(m)}_i{\bf 1}(x^1,x^2)1(x^1\!\not =\!x^2)\mu^1\!\times\!\mu^2(dx^1 dx^2).
\end{align}
and 
\begin{align}\label{cor3-e2}
\limsup\limits_{n\to\infty} \frac1{n} \ln \| \P_{\mu^1}(n,\cdot) - \P_{\mu^2}(n,\cdot)\|_{TV}
\!\le\! 
\liminf\limits_{m\to\infty} \limsup\limits_{n\to\infty} \frac1{nm} \ln \| \prod_{i=0}^{n-1} V^{(m)}_i \|.
\end{align}
The limit in the right hand side of (\ref{cor3-e2}) is uniform with respect to the initial measures $\mu^1, \mu^2$.
\end{corollary}
{\bf Proof} of (\ref{cor3-e1}) and (\ref{cor3-e2})  follows straightforwardly from (\ref{thm3-e1}). The right hand side in (\ref{cor3-e2}), clearly, does not depend on the initial measures. \hfill $\square$

\begin{remark}
In all bounds $ V^{(m)}$ may be replaced by $\hat V^{(m)}$; the values of the right hand sides will not change. Of course, the state space $ S^2$ is to be replaced by $\hat S^2$, and the product measure should be transformed respectively.
\end{remark}

\begin{remark}
Here is a remark on the famous Doob's estimate under the ``Doeblin -- Doob's condition (D)'' (\cite[chapter 5.5]{Doob}). Let us recall this  result in a slightly simplified version where only one ergodic subspace is allowed\footnote{In the original Doob's result there might be several or even infinite (countable) number of ergodic sets, so that the chain which starts at any state, eventually hits  one of them, after which the process can never escape from this ergodic set. We only show a simplified version with one ergodic set for comparison.}. The model is a homogeneous Markov chain in a general state space, with a transition kernel $Q(x,dx')$. Let the chain be acyclic and irreducible. Notation $Q_m(x,dx')$ stands for the $m$-step transition kernel, the same as $\mathsf P_x(m,dx')$.
\end{remark}

\begin{theorem}[J.L. Doob \cite{Doob}, theorem 5.5.6]\label{thm4}
If on the state space  $S$, there exist a probability measure $\lambda$, a natural number $m$, and $ \varepsilon >0$  such that from $\lambda(A)\le \varepsilon$ it follows $\sup_{x\in S} Q_m(x,A)\le 1- \varepsilon$, then there exist $C,c>0$ such that  
\begin{equation}\label{thm4-e1}
\sup_{x\in S} \|\mu^x_n - \mu_{inv}\|_{TV} \le C\exp(-cn), \quad n\ge 0.
\end{equation}
\end{theorem}  

\begin{remark}
The condition of the theorem \ref{thm4} excludes singularity of the $m$-step transition kernels for different starting states, say, $x_1$ and $x_2$. The result is exponential rate of convergence, and it may look like it could be compared to our theorems \ref{thm1}--\ref{thm3}, because the condition seems more general, while the rate is also exponential. However, the issue is that neither of the two constants $C$ and $c$ in (\ref{thm4-e1}) may be evaluated via $\varepsilon >0$ and $\lambda$ from the standing condition. In particular, $c>0$ may be arbitrarily close to zero, which was indicated by Doob himself. So, in fact, his result and the theorems \ref{thm1} -- \ref{thm3} of this paper are hardly comparable. Likely, the issue is because the measure $\lambda$ in the standing condition is not specified in the theorem \ref{thm4}, it is too abstract. On the contrary, in the theorems \ref{thm1} -- \ref{thm3} the transition kernels disclose the non-singularity in their own terms, which results in efficient bounds. 

\end{remark}

\section{Conclusive remarks}
For  Markov chains (MC) the new spectral approach for evaluating the rate of convergence in total variation uses a so called Markovian coupling. Namely, certain sub-stochastic matrices or operators, \(V\) and \(V^{(m)}\) have been introduced with the help of this coupling construction, and the rate of convergence under investigation in the homogeneous case has been estimated by their spectral radii \(r(V), r(V^{(m)})\). 

The point is that the ideal estimate for the convergence rate in total variation $\|\mu^x_n-\mu_\infty\|_{TV}$ is this total variation value itself, by definition. However, it is not a reasonable answer: we need some more simple tool for an efficiant estimate. Such a tool exists, it is Gantmacher's method based on the second by modulus eigenvalue of the transition matrix, see the proposition \ref{pro2}. It can be applied in the finite state space case. Another popular approach uses Markov -- Dobrushin ergodic coefficient, which is usually not optimal. What we claim is that ($1^\circ$) iterative Markov -- Dobrushin's ergodic coefficient may be as close to the best possible estimate as needed, by choosing an appropriate number of iterations, see the theorem \ref{cor}; and ($2^\circ$) there is even a better tool based on the Markovian coupling construction and on the spectral radious of the related sub-Markov kernel, \(r(V^{(m)})\), see the theorem \ref{thm2}. This method works well in any state space, not necessarily finite, and its advantage is, apparently, the existence of well-established procedures for evaluating spectral radius, see, for example, \cite[chapter 3, \S 16]{KLS}. Finally, the analogue of this latter method is developed in the non-homogeneous situation, see theorem \ref{thm3} and its corollary. 
The problem under consideration is also closely related to mixing coefficients: strong mixing, Kolmogorov's complete regularity, and uniformly strong mixing. In this respect, recall an exposition on mixing techniques for continuous time Markov processes in \cite{Ver-mix-note}.

\section*{Acknowledgments}
This research was carried out within the state assignment of Ministry of Science and Higher Education of the Russian Federation for IITP RAS. The author also thanks the anonymous referee whose remarks and comments helped make the text significantly better.

\end{document}